\documentclass{amsart}
\usepackage{amssymb}
\title{Multi-dimensional Weyl Modules\\ and Symmetric Functions}
\author{B. Feigin, S. Loktev}
\address{BF: Landau institute for Theoretical Physics, Chernogolovka
142432, Russia} \email{feigin@feigin.mccme.ru}
\address{SL:  Institute for Theoretical and Experimental Physics,
B. Cheremushkinskaja, 25, Moscow 117259, Russia}
\address{Independent University of Moscow, B. Vlasievsky per., 11,
Moscow 121002, Russia}
\email{loktev@mccme.ru}
\date{}
\newtheorem{dfn}{Definition}
\newtheorem{prp}{Proposition}
\newtheorem{lmm}{Lemma}
\newtheorem{thm}{Theorem}
\newtheorem{cnj}{Conjecture}
\newtheorem{crl}{Corollary}
\newcommand{\bin}[2]{\left[{#1}\atop{#2}\right]}
\begin{document}
\def\theenumi{\roman{enumi}}
\def\labelenumi{(\theenumi)}

\begin{abstract}
The Weyl modules in the sense of V.~Chari and A.~Pressley
(\cite{CP}) over the current Lie algebra on an affine variety are studied.
We show that local Weyl modules are finite-dimensional and
generalize the tensor product decomposition theorem from \cite{CP}.

More explicit results are stated for currents on a non-singular affine variety
of dimension $d$ with coefficients in the Lie algebra
$sl_r$. The Weyl modules with highest weights proportional to the vector representation one are related to
the multi-dimensional analogs of harmonic functions. The dimensions of such local Weyl modules
are calculated in the following cases.

For $d=1$ we show that the dimensions are equal to powers of $r$. 
For $d=2$ we show that the dimensions are given by products of the higher Catalan numbers 
(the usual Catalan numbers for $r=2$). We finally formulate a conjecture for an arbitrary $d$
and $r=2$.
\end{abstract}

\maketitle

\def \CC {{\mathbb C}}
\def \ZZ {{\mathbb Z}}
\def \g {{\mathfrak g}}
\def \h {{\mathfrak h}}
\def \n {{\mathfrak n}}
\def \b {{\mathfrak b}}
\def \ad {{\rm ad}\,}
\def \Z {{\mathcal Z}}

\section*{Introduction}

Let $M$ be a complex affine variety, let $\g$ be a semisimple Lie
algebra. By
$\g^M$ denote the Lie algebra of currents on $M$. For any point $\Z \in M$
we have the evaluation homomorphism $\g^M \to \g$ taking value at $\Z$. So
for any representation $\pi$ of $\g$ we can define the evaluation
representation $\pi(\Z)$ of $\g^M$.

\def \cn {{\mathcal N}}
\def \tn {\widetilde{\cn}}

Now let $\pi_1, \dots, \pi_n$ be irreducible representations of $\g$ and
let $\Z_1, \dots, \Z_n$ be points on $M$. Suppose that the points $\Z_i$
are
pairwise  distinct. Then $\pi_1(\Z_1) \otimes \dots \otimes \pi_n(\Z_n)$
is
an irreducible representation of $\g^M$. Fix the highest weight vectors $v_i \in
\pi_i$. Then the vector $v = v_1 \otimes \dots \otimes v_n$ is cyclic,
therefore
$$\pi_1(\Z_1) \otimes \dots \otimes \pi_n(\Z_n) \cong U(\g^M) / I(\pi_1,
\dots, \pi_n; \Z_1, \dots, \Z_n),$$
where $I(\pi_1, \dots, \pi_n; \Z_1, \dots, \Z_n)$ is the left ideal
annihilating $v$. So we have a family $\cn$ of ideals labeled by $M\times
\dots \times M$ without diagonals.

Let $\tn$ be the closure of $\cn$ in the space of all left ideals in
$U(\g^M)$. Then $\tn$ is an algebraic variety (in most cases singular).
Quotients of $U(\g^M)$ by the corresponding ideals form a bundle $\xi$ on
$\tn$ of dimension $\dim (\pi_1) \cdots \dim (\pi_n)$.

The case $M=\CC$, $\g = sl_2$ is completely studied in \cite{FF}. Here the limits of
ideals are unique and
the fiber over $(0,\dots,0)$ is isomorphic to the {\em fusion product}
introduced in \cite{FL}. Now suppose that all the
representations are 2-dimensional.
Then $\tn$ is just the $n$-th symmetric power of $\CC$ and
the fusion product is universal in
the following sense. It is the maximal finite-dimensional module generated
by the vector $v$ such that
$$(e\otimes P)\, v = 0, \qquad (h\otimes P)\, v = n P(0) \cdot v,$$
where $P$ is a polynomial, and $e$, $h$, $f$ constitute the standard basis
in $sl_2$. Following \cite{CP} we call such representations {\em Weyl
modules}.

We expect that for $M = \CC^d$, $\g = sl_2$ and 2-dimensional
representations the variety $\tn$ is a general locus in the Hilbert scheme of
$n$-points on $\CC^d$. So there is the canonical projection $\theta : \tn
\to S^n(\CC^d)$. By $\tn_0$ denote the special fiber $\theta^{-1}
(0,\dots, 0)$. Then we expect that the space of sections of the
$2^d$-dimensional bundle $\xi$ on $\tn_0$ has the same universality
property.

Note that both facts can be derived for $d=2$ from the recent results of
M.~Haiman (\cite{H}). It gives a very powerful instrument to study
fusion and Weyl modules over two-dimensional currents.

In this paper we concentrate on the Weyl module side of this picture. In
Section~1 we extend the notion of Weyl modules.
Namely we define the {\em global} Weyl module corresponding to a
commutative algebra with unit, a semisimple Lie algebra and a weight, as well as
the {\em local} Weyl module corresponding to the data above and a set of
points
on ${\rm Spec}\ A$. We show that local Weyl modules are finite-dimensional
and generalize the tensor product theorem (\cite{CP}) for them.
In Section~2 we generalize the symmetric power construction proposed
in \cite{CP}. In Section~3 we apply the results of \cite{H} to 
investigation of
Weyl modules over two-dimensional currents 
and state a conjecture for the general case.

\subsection*{Acknowledgments}
This consideration is inspired by discussions with Vyjayanthi Chari at UC
Riverside and MSRI. SL thanks V.~Chari for the hospitality at UC
Riverside. We are also grateful to
V.~V.~Dotsenko, A.~N.~Kirillov and I.~N.~Nikokoshev for useful and stimulating discussions.

BF is partially supported by the grants
RFBR-02-01-01015, RFBR-01-01-00906 and INTAS-00-00055.

SL is partially supported by the grants RFBR-02-01-01015,
RFBR-01-01-00546 and CRDF RM1-2545-MO-03.

\section{Some generalities}

\subsection{Notation and integrability}

Let $\g$ be a semisimple Lie algebra of rank $r$. Choose a Cartan subalgebra
$\h$ and
Borel subalgebras $\b = \b_+$ and $\b_-$. Denote their nilpotent radicals
by  $\n_+$ and $\n_-$ respectively.

By $Q\subset \h^*$ denote the root lattice, by $P\subset \h^*$ denote the
weight lattice. 
Let $\alpha_1, \dots, \alpha_r \in Q$ be the simple roots,
let $\omega_1, \dots, \omega_r \in P$ be the simple weights.
Denote the correspondent Chevalley generators by
$e_1, \dots, e_r \in \n_+$, $h_1, \dots, h_r \in \h$ and
$f_1, \dots, f_r \in \n_-$.
Let $Q_+ \subset Q$ be the cone of positive roots, that is the set of 
linear combination of simple roots $\alpha_i$ with non-negative 
integer coefficients.

\begin{dfn}
We say that a $\g$-module $V$ is {\em integrable} if there exists a finite subset $D \subset P$
such that
$$V = \bigoplus_{\mu \in D} V^\mu,$$
where any $h \in \h$ acts on $V^\mu$ by the scalar $\mu(h)$.
\end{dfn}

Note that for our purposes we do not require $V^\mu$ to be
finite-dimensional.  Similar to the case of finite-dimensional weight spaces (see \cite{K})
it actually implies that 
$V$ is a possibly infinite sum of finite-dimensional representations. 

\begin{prp}\label{weyl}
There is the action of the Weyl group on the weights of integrable representations.
\end{prp}

\begin{proof}
It is enough (see \cite{K}) to introduce the action of the operators
$$s_i = \exp(f_i) \exp(-e_i) \exp(f_i)$$
corresponding to simple reflections.
These operators are well-defined on an integrable representation because
the action of $e_i$ and $f_i$ is nilpotent. So similar to the case of finite-dimensional 
representations it defines the action of the Weyl group on the weights.
\end{proof}

Also we consider the following class of representation.

\begin{dfn}
Given a weight $\lambda \in \h^*$ we say that a $\g$-module $V$ is {\em $\lambda$-bounded} if
$$V = \bigoplus_{\mu \in \lambda - Q_+} V^\mu,$$
where any $h \in \h$ acts on $V^\mu$ by the scalar $\mu(h)$.
\end{dfn}

This definition is related to the definition of the category ${\mathcal O}$ (see \cite{K}). 
Namely any representation from the category ${\mathcal O}$ is by definition a direct sum of
$\lambda^i$-bounded modules for some $\lambda^1, \dots, \lambda^N \in \h^*$
with finite-dimensional weight spaces.

\begin{prp}
Any $\lambda$-bounded $\g$-module has the maximal integrable quotient.
\end{prp}

\begin{proof}
Suppose that a module $V = \bigoplus V^\mu$ is $\lambda$-bounded. By $w_0 \in W$ denote
the highest length element, so $w_0 (Q_+) = - Q_+$. 
Let $V'$ be the submodule of $V$ generated by $V^\mu$
such that $\mu \not\in w_0(\lambda) + Q_+$. Introduce the module $V^{int}
= V / V'$. As the intersection of $\lambda - Q_+$ and $w_0(\lambda) + Q_+$ is finite, the module
$V^{int}$ is integrable.

Now let $U$ be an integrable quotient of $V$. Let us show that $U$ is
indeed a quotient of $V^{int}$.

As $U$ is a quotient of $V$ we have
$U^\mu =0$ if $\mu \not\in  \lambda - Q_+$. Then 
Proposition~\ref{weyl}
implies that
$U^\mu =0$ if $\mu \not\in w_0(\lambda -Q_+) = w_0(\lambda) + Q_+$. 
Therefore the image of $V'$ in $U$
is zero, so $U$ is
a quotient of $V^{int}$.
\end{proof}

Let us say that a module with a certain property is
{\em maximal} if any other module with this property is a quotient of
the given one.

\subsection{Weyl modules}

\def \ind {{\rm Ind\,}}
\def \hom {{\rm Hom}}

Let $A$ be a commutative algebra with unit. Generalizing \cite{CP} we
define a family of representations of the Lie algebra $\g \otimes A$.

\begin{dfn} For $\lambda \in \h^*$  define the
{\em global Weyl} module
$W_A(\lambda)$ as the
maximal $\g\otimes 1$-integrable module generated by the vector
$v_\lambda$ such
that
\begin{equation}\label{hwc}
(\n_+ \otimes P) \, v_\lambda = 0, \qquad (h\otimes 1) \, v_\lambda =
\lambda(h) \cdot v_\lambda \quad \mbox{for} \ h \in \h, \ P \in A.
\end{equation}
\end{dfn}

In particular if $A = \CC[x,x^{-1}]$ then $W_A (\lambda) \cong
W(\lambda)$ in the
terms of \cite{CP}.

\begin{prp}\label{welld1}
The global Weyl module $W_A(\lambda)$ exists.
\end{prp}

\begin{proof}
Let $\b_A \subset \g\otimes A$ be the subalgebra spanned by $\n_+ \otimes
A$ and $\h \otimes 1$. Then the induced module
$$M_A(\lambda) = \ind_{\b_A}^{\g\otimes A} \, v_\lambda$$
is the maximal module generated by $v_\lambda$.

As $M_A(\lambda)$ is
generated by
$v_\lambda$ under the action of $U(\b_- \otimes A)$, the module
$M_A(\lambda)$ is $\lambda$-bounded. So
we can define
$W_A(\lambda)$ as the maximal
integrable quotient of $M_A(\lambda)$.

Concerning the universality property, any integrable module generated by $v_\lambda$ is an 
integrable quotient of $M_A(\lambda)$ and therefore a quotient of the constructed module
$W_A(\lambda)$.
\end{proof}

\begin{lmm}
We have $W_A(\lambda) \ne 0$ if and only if $\lambda$ is dominant.
\end{lmm}

\begin{proof}
Note that $v_\lambda$ can generate an integrable module only for
a dominant $\lambda$. 

On the other hand, suppose $\lambda$ is dominant. Let $\pi_\lambda$ be the irreducible 
representation of $\g$ with highest weight $\lambda$, let
$v$ be its highest weight vector. Then $\pi_\lambda \otimes A$ is 
an integrable module over $\g \otimes A$ and $v\otimes 1$ satisfies the conditions~\eqref{hwc}.
So $W_A(\lambda)$ has a non-trivial quotient $\pi_\lambda \otimes A$, 
therefore $W_A(\lambda) \ne 0$.
\end{proof}

\def \sz {{\{\Z\}_\lambda}}

Now regard $A$ as the algebra of functions on an affine
variety $M$ (possibly singular). Actually points on $M$ are just
ideals in
$A$ of codimension one.

Now suppose that $\lambda$ is dominant. Then the coordinates
$\lambda_i =\lambda(h_i)$ are non-negative integers.
A set of points $\Z^{(i)}_j \in M$ (not necessarily distinct)
indexed as follows:
$$\sz = \left(\{\Z^{(1)}_j\}_{j=1\dots \lambda_1}, \dots,
\{\Z^{(r)}_j\}_{j=1\dots
\lambda_r}\right)
$$
we call a {\em $\lambda$-multiset}.

\begin{dfn}
Let $\lambda = (\lambda_1, \dots, \lambda_r)\in \h^*$ be a
dominant weight, let $\sz$ be a $\lambda$-multiset. Define the
{\em local Weyl} module
$W_M (\sz)$ as the
maximal $\g\otimes 1$-integrable module generated by the vector $v_\sz$
such that
\begin{equation}\label{hwl}
(\n_+ \otimes P) \, v_\sz = 0, \qquad (h_i\otimes P) \, v_\sz =
\sum_{j=1}^{\lambda_i} P\left(\Z^{(i)}_j\right) \cdot v_\sz \quad
\mbox{for} \ i=1\dots r, \ P \in A.
\end{equation}
\end{dfn}

In particular if $M = \CC^*$, $A = \CC[x,x^{-1}]$ then $W_M (\sz) \cong
W(\pi)$ in the
terms of \cite{CP} for $\pi_i(u) = \prod_j \left(1-u/\Z^{(i)}_j\right)$.
In general case this definition is motivated by Theorem~\ref{univ} from
Section~2.

\begin{prp}
The local Weyl module $W_M(\sz)$ exists and 
can be obtained as a quotient of the
global Weyl module $W_A(\lambda)$.
\end{prp}

\begin{proof}
Similarly to Proposition~\ref{welld1} one can define $W_M(\sz)$ as 
the maximal integrable quotient of
the induced module
$$M_M(\sz) = \ind_{\b\otimes A}^{\g\otimes A} \, v_\sz.$$

As $\b_A$ acts on $v_\sz$ in the same way as on $v_\lambda$, the
universality property implies that $W_M(\sz)$ is a
quotient of $W_A(\lambda)$.
\end{proof}

\begin{prp}
We have $W_M(\sz) \ne 0$.
\end{prp}

\begin{proof}
With any $\Z \in \sz$ we associate the weight $\lambda_\Z$ such that $h_i(\lambda_\Z)$ is the 
multiplicity of $\Z$ in $\{\Z^{(i)}_j\}_{j=1\dots\lambda_i}$.

For $\Z \in M$ and a representation $\pi$ of $\g$ let $\pi(\Z)$ be the evaluation representation
of $\g \otimes A$ at $\Z$. Namely $\pi(\Z) = \pi \otimes (A/ I(\Z))$, where $I(\Z)$ is 
the maximal ideal corresponding to the point $\Z$.

By $\pi_\mu$ denote the irreducible representation with highest weight $\mu$.
Then the representation 
$$\bigotimes_{\Z \in \sz} \pi_{\lambda_\Z}(\Z)$$ 
is integrable and the product of highest weight vectors satisfies~\eqref{hwl}. So $W_M(\sz)$ has a
non-trivial quotient, therefore  $W_M(\sz) \ne 0$.
\end{proof}

At last note that the action of the Weyl group restricts weights of the Weyl modules as follows.

\begin{lmm}\label{pgr}
We have $W_A(\lambda)^\mu \ne 0$ only for 
$\mu \in \bigcap\limits_{w \in W} w(\lambda -Q_+)$, and the same statement for $W_M(\sz)$.
\end{lmm}

\begin{proof}
We just combine the observation that $W_A(\lambda)$ and  $W_M(\sz)$ are $\lambda$-bounded
and the action of the Weyl group on their weights.
\end{proof}

\subsection{The Local Weyl modules are finite-dimensional}

\begin{lmm}\label{igr}
We have $W_A(\lambda)^{\lambda-N\alpha_i} \ne 0$ only for $0\le N \le \lambda_i$, 
and the same statement for $W_M(\sz)$.
\end{lmm}

\begin{proof}
It follows from Lemma~\ref{pgr}. For $N<0$ we have $\lambda-N\alpha_i \not\in \lambda -Q_+$,
for $N>\lambda_i$ we have $\lambda-N\alpha_i \not\in s_i(\lambda -Q_+)$, where $s_i$ is the simple
reflection.
\end{proof}

Let $I(\sz)$ be the ideal in $A$ of functions with zeroes at all
$\Z^{(i)}_j$. In order to get multiple zeroes we consider the powers
$I^s (\sz)$.

\begin{prp}\label{anr}
Elements of
$f_i \otimes I^{\lambda_i}(\sz)$ act on $v_\sz$ by
zero.
\end{prp}

\begin{proof}
The ideal $I^{\lambda_i}(\sz)$ is spanned by decomposable elements, that is
by products of
$\lambda_i$ elements of $I(\sz)$. So it is enough to prove that if
$P = P_1 \cdots
P_{\lambda_i}$ for $P_i \in I(\sz)$ then $(f_i \otimes P) \, v_\sz =0$.

From Lemma~\ref{igr} it follows that $(f_i\otimes 1)^{\lambda_i+1} v_\sz
=0$.
Then using Proposition~\ref{martini} (see Appendix) we have 
$$0 = (e_i \otimes P_1)\cdots (e_i\otimes P_{\lambda_i}) \cdot
(f_i\otimes 1)^{\lambda_i+1} v_\sz = 
C\cdot (f_i \otimes P) \, v_\sz,$$
where the coefficient $C \ne 0$.
\end{proof}

\begin{prp}\label{anid}
For a large enough integer $N$ the ideal $\g \otimes I^N(\sz)$ acts on
$W_M(\sz)$ by zero.
\end{prp}

\begin{proof}
Note that the Lie subalgebra $\n_-$ is generated by $f_i$, $i=1\dots
r$. So $\n_-$ is spanned by a finite number of elements $x_{(i_1, \dots,
i_s)} = \ad(f_{i_1})
\cdots  \ad(f_{i_{s-1}})\, f_{i_s}$. Let us take
$$N \ge \max_{{(i_1, \dots, i_s)}\atop{x_{(i_1, \dots, i_s)}\ne 0}}
(\lambda_{i_1} +
\dots + \lambda_{i_s}).$$
Then for a decomposable $P \in I^N(\sz)$
and $(i_1, \dots, i_s)$ such that $x_{(i_1,
\dots, i_s)}\ne 0$ we have $P = P_1 \cdots  P_s$, where
$P_j \in  I^{\lambda_{i_j}}(\sz)$. From Proposition~\ref{anr} we have
$(f_{i_j}\otimes P_j)\, v\sz$ =0, so
$$(x_{(i_1, \dots,
i_s)} \otimes P) \, v_\sz = \left(\ad(f_{i_1}\otimes P_1)
\cdots  \ad(f_{i_{s-1}}\otimes P_{s-1}) (f_{i_s}\otimes P_s)\right)
\, v_\sz =0.$$
It means that $\n_- \otimes I^N(\sz)$ acts on $v\sz$ by zero.

As $\b_+ \otimes I(\sz)$ acts on $v_\sz$ by zero and $\g$ = $\b_+ \oplus
\n_-$, we have $(\g \otimes
I^N(\sz))\, v_\sz =0$. Finally note that $\g
\otimes
I^N(\sz)$ is an ideal in the Lie algebra $\g\otimes A$, so we have the
statement of the
proposition.
\end{proof}

\begin{thm}
The local Weyl module $W_M(\sz)$ is finite-dimensional.
\end{thm}

\begin{proof}
Recall that  $W_M(\sz)$ is generated by $v_\sz$ under the action of $\n_-
\otimes A$, in other words we have   $W_M(\sz) = U(\n_- \otimes A) \, v_\sz$. From
Proposition~\ref{anid} we have that $W_M(\sz) = U\left(\n_- \otimes
A/I^N(\sz)\right) \, v_\sz$ for a certain $N$.

\def \rv {{\rho^\vee}}

Now let us deduce from integrability that there exists an integer $n$ such that
$x_1 \cdots x_s \, v_\sz =0$ in $W_M(\sz)$ for $s > n$ and any  $x_i \in \n_-
\otimes A$. 

Let $\rv$ be the element of $\h$ such that $\rv(\alpha_i) = 1$ for
$i=1\dots r$. As $W_M(\sz)$ is integrable, the set of values of $\rv$ on weights of
$W_M(\sz)$ is finite.
Suppose that $x_1 \cdots x_s \, v_\sz \ne 0$, let
$\mu $ be its weight. Then $\rv(\mu) < \rv(\lambda) - s$.
So $s$ is not greater than the maximal difference between $\rv(\lambda)$
and $\rv(\mu)$, where $\mu$ is a weight of $W_M(\sz)$.

Now we have
$$W_M(\sz) = U^{\le n}\left(\n_- \otimes
A/I^N(\sz)\right) \, v_\sz.$$
As the Lie algebra $\n_- \otimes A/I^N(\sz)$ is finite-dimensional, this
space is also finite-dimensional.
\end{proof}

\subsection{Tensor product decomposition}

For two multisets $\sz$ and $\{\Z'\}_{\lambda'}$ define the union multiset
$$\{\Z \cup \Z'\}_{\lambda+\lambda'} = \left(\{\Z^{(1)}_j\}\cup
\{(\Z')^{(1)}_j\}, \dots,
\{\Z^{(r)}_j\}\cup \{(\Z')^{(r)}_j\}\right).$$

\begin{lmm}\label{t1}
We have an isomorphism of 1-dimensional representations of $\b \otimes A$
$$v_{\{\Z \cup \Z'\}_{\lambda+\lambda'}} \cong v_\sz \bigotimes
v_{\{\Z'\}_{\lambda'}}.$$
\end{lmm}

\begin{proof}
It follows
from the definition~\eqref{hwl}.
\end{proof}

\begin{prp}\label{s1}  Let $A_{N\sz} = A/ I^N_\sz$.
Suppose that $\Z^{(i)}_j \ne (\Z')^{(i')}_{j'}$ for all $i$, $j$, $i'$,
$j'$.
Then we have an isomorphism of algebras
$$A_{N\{\Z \cup \Z'\}_{\lambda+\lambda'}} \cong A_{N\sz} \bigoplus
A_{N\{\Z'\}_{\lambda'}}.$$
\end{prp}

\begin{proof}
From the Hilbert Nullstellensatz we have that $I^N_{\{\Z \cup
\Z'\}_{\lambda+\lambda'}} = I^N_\sz \cap I^N_{\{\Z'\}_{\lambda'}}$ and
that $I^N_\sz + I^N_{\{\Z'\}_{\lambda'}} = A$. Then the direct sum of
the natural projections
$$A \to  A_{N\sz} \bigoplus A_{N\{\Z'\}_{\lambda'}}$$
is surjective and its kernel is equal to $I^N_{\{\Z \cup
\Z'\}_{\lambda+\lambda'}}$.
\end{proof}

\begin{thm}\label{tcmp}
Suppose that $\Z^{(i)}_j \ne (\Z')^{(i')}_{j'}$ for all $i$, $j$, $i'$,
$j'$.
Then we have an isomorphism of $\g \otimes A$-modules
$$W_M(\{\Z \cup \Z'\}_{\lambda+\lambda'}) \cong W_M(\sz) \bigotimes
W_M(\{\Z'\}_{\lambda'}).$$
\end{thm}

\begin{proof}
By Lemma~\ref{t1} we have that
the product $W_M(\sz) \bigotimes
W_M(\{\Z'\}_{\lambda'})$ is a quotient of $W_M(\{\Z \cup
\Z'\}_{\lambda+\lambda'})$. Let us show the equality of their dimensions. 

Right now we know that the dimension of $W_M(\{\Z \cup
\Z'\}_{\lambda+\lambda'})$ is not less than 
the dimension of the tensor product. Let us obtain
the opposite inequality.

Recall the notation $A_{N\sz} = A/ I^N_\sz$ and
consider the
corresponding induced module
$$M^N_M(\sz) = \ind_{\b\otimes A_{N\sz}}^{\g\otimes A_{N\sz}} \, v_\sz.$$
From Proposition~\ref{anid} we know that $W_M(\sz)$ is 
a quotient of $M^N_M(\sz)$ for a sufficiently large $N$. Moreover as
$M^N_M(\sz)$
is a quotient of $M_M(\sz)$, we have that $W_M(\sz)$ is its maximal
integrable quotient.

Lemma~\ref{t1} with Proposition~\ref{s1}
implies that for any positive integer
$N$ we have
$$M^N_M(\{\Z \cup \Z'\}_{\lambda+\lambda'}) \cong M^N_M(\sz) \boxtimes
M^N_M(\{\Z'\}_{\lambda'}).$$
Here the right hand side is the external tensor product of a $\g \otimes A_{N\sz}$-module and a
 $\g \otimes A_{N\{\Z'\}_{\lambda'}}$-module. 

\def \id {{\rm Id}}

    Now take a sufficiently large $N$ so three our $W_M$-s are quotients of 
the corresponding $M^N_M$-s. 
Then the module $W_M(\{\Z \cup\Z'\}_{\lambda+\lambda'})$ 
is a quotient of the tensor product $M^N_M(\sz) \boxtimes  
M^N_M(\{\Z'\}_{\lambda'})$.

Concerning the weight subspaces, let us show that the image of 
$M^N_M(\sz)^\mu \boxtimes M^N_M(\{\Z'\}_{\lambda'})^{\mu'}$ in
the quotient $W_M(\{\Z \cup \Z'\}_{\lambda+\lambda'})$ 
is zero except  a finite number of $\mu$ and $\mu'$.

Let $D$ be the set of weights of the module $W_M(\{\Z \cup
\Z'\}_{\lambda+\lambda'})$. As this module is integrable, $D$ is finite.
Then the sets
$$D_1 = (\lambda - Q_+) \bigcap (-\lambda'+D+Q_+), \qquad 
D_2 = (\lambda' - Q_+) \bigcap (-\lambda+D+Q_+)$$
are also finite.
As the module $M^N_M(\sz)$ is $\lambda$-bounded and $M^N_M(\{\Z'\}_{\lambda'})$ is
$\lambda'$-bounded, the image of $M^N_M(\sz)^\mu \boxtimes M^N_M(\{\Z'\}_{\lambda'})^{\mu'}$ in $W_M(\{\Z \cup
\Z'\}_{\lambda+\lambda'})$ is zero unless $\lambda -\mu \in Q_+$, $\lambda' -\mu' \in Q_+$
and, of course, $\mu+\mu' \in D$. So it is not zero only if $\mu \in D_1$ and $\mu' \in D_2$.

Let $M^N_M(\sz)'$ be the submodule of $M^N_M(\sz)$ generated by 
weight subspaces $M^N_M(\sz)^\mu$ 
with $\mu \not\in D_1$. Consider the quotient $W^N_M(\sz) = M^N_M(\sz)/M^N_M(\sz)'$.
Let $W^N_M(\{\Z'\}_{\lambda'})$ be the similar quotient of 
$M^N_M(\{\Z'\}_{\lambda'})$
with weights in  $D_2$. 
Then we know that the 
module
$W_M(\{\Z \cup
\Z'\}_{\lambda+\lambda'})$ is a quotient of $W^N_M(\sz) \boxtimes W^N_M(\{\Z'\}_{\lambda'})$.

Because the sets $D_1$ and $D_2$ are finite, 
the modules $W^N_M(\sz)$ and $W^N_M(\{\Z'\}_{\lambda'})$
 are integrable. Then the universality property implies that 
the module  $W_M(\{\Z \cup
\Z'\}_{\lambda+\lambda'})$ is indeed a quotient of $W_M(\sz) \boxtimes
W_M(\{\Z'\}_{\lambda'})$. So 
$$\dim \left(W_M(\sz) \bigotimes
W_M(\{\Z'\}_{\lambda'})\right) \ge \dim W_M(\{\Z \cup
\Z'\}_{\lambda+\lambda'}),$$
that completes the proof.
\end{proof}

\begin{crl}\label{elw}
Suppose that $M$ is non-singular of dimension $d$. Then any local Weyl
module is isomorphic as a vector space
to the tensor product of modules $W_{\CC^d}(\{0\}_{\lambda^i})$ for some
dominant $\lambda^1, \dots, \lambda^s$.
\end{crl}

\begin{proof}
By Theorem~\ref{tcmp} it is enough to describe $W_M(\sz)$ when $\sz$
contains only one point $Z \in M$. In this case $I(\sz)$ is the ideal
$I_Z$ of
functions with zero at $Z$. By Proposition~\ref{anid} we know that starting from a certain $N$
$W_M(\sz)$ depends only on $A/ I^N(\sz)$. And the algebra $A/I^N_Z$ is the same for
any non-singular point on any affine variety of given dimension.
\end{proof}

To proceed with a singular variety $M$ we need also to describe local Weyl modules at
singular points.

\section{The symmetric power construction}

Here we generalize the results of Section~6 in \cite{CP}. Our proofs are based 
on the following lemma.

\begin{lmm}\label{ordspan}
Let $X$ be a vector space with a basis $v_i$, $i\in I$. Suppose that the set 
$I$ is partially ordered
such that any chain $i_1 > i_2 > i_3 > \dots$ is finite.
\begin{enumerate}
\item Let $\{u_i\}$, $i\in I$, be a set of vectors such that each $u_j$ is the sum of $v_j$ 
and a linear combination of $v_i$ with $i <j$. Then the set $\{u_i\}$ is a basis of $X$.
\item Let $Y\subset X$ be a subspace.
Suppose we have a subset $J \subset I$ such that for any
$j \in J$ the vector $v_j$ plus a linear combination of $v_i$ with $i < j$ belongs to $Y$.  
Then the images of vectors $v_i$ for $i \not\in J$ span $X/Y$.
\end{enumerate}
\end{lmm}

\begin{proof}
To prove (i) it is enough to show 
that each $v_j$ is a linear combination of $u_i$ as well.
The vector $v_j$ is equal to $u_j$ plus (maybe) a linear combination 
of some $v_{j'_a}$ with $j'_a < j$.
Similarly each   $v_{j'_a}$  is equal to $u_{j'_a}$ plus (maybe) a linear combination 
of some $v_{j''_{ab}}$ with $j''_{ab} < j'_a$ and so on. 
As any decreasing chain is finite, this procedure is also finite, and the result is
a linear combination of $u_i$.

To show (ii) let us take $u_i = v_i$ for $i \not\in J$ and 
$u_i = v_i + \dots \in Y$ for $i \in J$. Then by (i) the set $\{u_i\}$ is a basis, so any vector
from $X$ is a linear combination of some $v_i$ for $i \not\in J$ and a vector from $Y$.
\end{proof}

\subsection{The $sl_2$ case} For $\g=sl_2$ let us write $e$, $h$, $f$
instead of $e_1$, $h_1$, $f_1$.

\def \natm {S^n \left(L_A(V)\right)}

Let $V$ be the vector (2-dimensional) representation of $sl_2$.
Then the vector space $L_A(V) = V \otimes A$ inherits the action of
$sl_2\otimes A$. Following \cite{CP} we consider the symmetric
powers $\natm$.

Choose a basis $\{ u, v\}$ in $V$ such that $hv =v$, $hu = -u$. Then we
can decompose $\natm$ into eigenspaces of $h\otimes 1$.

\begin{prp}\label{decmp}
We have 
$$
\natm \cong \bigoplus_{i=0}^n \natm^{n-2i}, 
\qquad \mbox{where}\ \ 
\natm^{n-2i} \cong S^i(A) \bigotimes S^{n-i}(A) \quad \mbox{for} \ \  i=0\dots n
$$
and $h\otimes 1$ acts on $\natm^s$ by the scalar $s$. 
\end{prp}

\begin{proof}
We have
$$\natm \cong S^n\left((u \otimes A) \bigoplus (v \otimes A)\right) \cong
\bigoplus_{i=0}^n S^i(u\otimes A) \bigotimes S^{n-i}(v\otimes A) \cong
\bigoplus_{i=0}^n S^i(A)\bigotimes S^{n-i}(A).$$
\end{proof}

Let
$$v^{(n)} = (v\otimes 1) \bigotimes  \dots \bigotimes (v\otimes 1).$$
This vector is cyclic in the module $\natm$. Now
choose a basis $\{P_i\}$ in $A$ such that $P_1 =1$.
To avoid cumbersome notation let us write $gP_i$ instead of $g \otimes P_i$
and $g$ instead of $g\otimes 1$ for $g \in sl_2$.

\begin{prp}\label{bas}
Monomials
$$fP_{j_1} \cdots  fP_{j_i} \cdot hP_{j_{i+1}} \cdot \dots
\cdot hP_{j_n} \, v^{(n)}, \quad j_1 \ge \dots \ge j_i, \quad j_{i+1}
\ge \dots \ge j_n, \quad i=0\dots n$$
form a basis in $\natm$.
\end{prp}

\begin{proof}
By definition the elements
$$s^i_{j_1, \dots, j_n} = {\rm Sym} \left((u\otimes P_{j_1}) \bigotimes 
\dots \bigotimes
(u\otimes P_{j_i}) \bigotimes (v
\otimes P_{j_{i+1}})\bigotimes \dots \bigotimes (v\otimes 
P_{j_n})\right)$$
for
$j_1 \ge \dots \ge j_i$, $j_{i+1} \ge \dots \ge j_n$
and  $i=0\dots n$
form a basis of $\natm$.

Introduce the partial order on the vectors $(j_1, \dots, j_n)$ 
that compares the number of $j_k$ 
equal to one for $k>i$:
$$(j_1, \dots, j_n) > (j_1', \dots, j_n')\quad \mbox{if} \quad
\left|\{ k > i | j_k =1\}\right| < \left|\{ k > i | j_k' =1\}\right|.$$
Clearly it satisfies the condition of Lemma~\ref{ordspan}.

Then we have
$$fP_{j_1} \cdots  fP_{j_i} \cdot hP_{j_{i+1}} \cdots
 hP_{j_n} \, v^{(n)} = n! \cdot s^i_{j_1, \dots, j_n} + \mbox{lower terms},$$
where lower terms are given by a linear combination of $s^i_{j'_1, \dots, j'_n}$ with
$(j_1', \dots, j_n') < (j_1, \dots, j_n)$.
Namely the main term appears when the factors in the left hand side are applied to
pairwise distinct factors of $v^{(n)}$, and lower terms 
appear when
some factors in the left hand side are
applied to the same factor of $v^{(n)}$.

By Lemma~\ref{ordspan}~(i) we have the statement of the proposition.
\end{proof}

\def \vn {{v_{n\omega}}}

\begin{prp}\label{span}
The global Weyl module $W_A(n\omega)$ is  spanned by  elements
\begin{equation}\label{mon1}
fP_{j_1} \cdots  fP_{j_i} \cdot
hP_{j_{i+1}} \cdots hP_{j_n} \, \vn, \quad j_1 \ge \dots \ge j_i, 
\quad j_{i+1}
\ge \dots \ge j_n, \qquad i=0\dots n.
\end{equation}
\end{prp}

\begin{proof}
The PBW Theorem and Lemma~\ref{igr} imply that $W_A(n\omega)$ is  spanned by  elements
\begin{equation}\label{mon2}
m^i_{j_1, \dots, j_s}=
fP_{j_1} \cdots fP_{j_i} \cdot hP_{j_{i+1}} \cdots
hP_{j_s} \, \vn,  \quad j_1 \ge \dots \ge j_i, 
\quad j_{i+1}
\ge \dots \ge j_s, \qquad i \le n.
\end{equation}
We need to show that this is still true under the
restriction $s\le n$. Then by setting $j_{s+1} = \dots = j_n = 1$ 
(that is $P_j=1$ for $j>s$) we obtain the statement of the proposition.

Now fix the grading, that is fix the integer $i$. 
Introduce the partial order 
$$(j_1, \dots, j_s) > (j_1', \dots, j_{s'}')\quad \mbox{if} \ s>s' \ \mbox{or}
\left( s=s' \ \mbox{and} \ \left|\{ k \le i | j_k =1\}\right| < \left|\{ k \le i | j_k' =1\}\right|
\right).$$
Clearly it satisfies the condition of Lemma~\ref{ordspan}.

To show that the monomials~\eqref{mon1} span $W_A(n\omega)$ let us just obtain
the monomials~\eqref{mon2} with $s>n$ as the highest terms of some relations in
$W_A(n\omega)$.

As $eP \cdot \vn =0$ for all $P \in A$, we have
$$eP \cdot m^i_{j_1, \dots, j_s} = \sum_{k=1}^i
fP_{j_1} \dots fP_{j_{k-1}} fP_{j_{k+1}} \dots fP_{j_i}
h(PP_{j_k}) hP_{j_{i+1}}  \dots hP_s \vn + \mbox{lower length terms},$$
where lower length terms are given by a linear combination of $m^{i'}_{j'_1, \dots, j'_{s'}}$
with $s'<s$. Therefore we have
\begin{equation}\label{rel}
eP_{j_{i+1}} \cdots eP_{j_s} \cdot f^{s-i} \cdot fP_{j_1} \cdots
fP_{j_i} \, \vn =
(s-i)! \, m^i_{j_1, \dots, j_s} + \mbox{lower terms},
\end{equation}
where lower terms are given by a linear combination of $m^{i'}_{j'_1, \dots, j'_{s'}}$ with
$(j_1', \dots, j_{s'}')< (j_1, \dots, j_s)$.

Now let $s>n$. Then from Lemma~\ref{igr} we have $f^{s-i} \cdot fP_{j_1}
\cdots fP_{j_i} \, \vn =0$, therefore the left hand side of~\eqref{rel}
is a relation in $W_A(n\omega)$
and $m^i_{j_1, \dots, j_s}$ is the highest term of this relation.

Let $X$ be the formal span of monomials~\eqref{mon2} and let $Y$ be the kernel of the natural map
$X \to W_A(n\omega)$.
Then Lemma~\ref{ordspan}~(ii) implies the statement of the proposition.
\end{proof}

\begin{thm}\label{iso}
We have the natural isomorphism of $sl_2 \otimes A$-modules
$$W_A(n\omega) \cong \natm.$$
\end{thm}

\begin{proof}
Note that $v^{(n)}$ satisfies~\eqref{hwc}, and that $\natm$ is generated by
 $v^{(n)}$.
Then from the universality property we have the surjective map $W_A(n\omega)
\to \natm$ sending $\vn$ to $v^{(n)}$. Propositions~\ref{bas} and
\ref{span} imply that this map is an isomorphism.
\end{proof}

\def \szn {{\{\Z\}_{n\omega}}}

Note that the action of $A$ on itself gives us the action of $S^n(A)$ on
$\natm$ commuting with the action of $sl_2 \otimes A$. For a multiset $\szn$
define the ideal
$$J_\szn = \left\{ R \in S^n(A) \left| \,
R\left(\Z^{(1)}_1, \dots, \Z^{(1)}_n\right) =
0\right.\right\}.$$

\begin{thm}\label{iso2}
For any $n\omega$-multiset $\szn$
we have the natural isomorphism of $sl_2 \otimes A$-modules
$$W_M(\szn) \cong \natm / J_\szn \cdot \natm.$$
\end{thm}

\begin{proof}
From Theorem~\ref{iso} and the universality property we have a surjective map $\natm \to
W_M(\szn)$ sending $v^{(n)}$ to $v_\szn$.

First let us show that the image of $J_\szn \cdot \natm$ under this map is
zero. As the action of $S^n(A)$ commutes with the action of $\g \otimes A$,
it is enough to show that the image of $J_\szn \cdot v^{(n)}$ is zero.

Note that the tensor algebra $T^*(A) \cong U(h\otimes A)$ acts on
$v^{(n)}$. It defines the surjective map
$$T^*(A) \to \natm^\lambda \cong S^n(A): \qquad
P_1\bigotimes \dots \bigotimes P_s \to hP_1 \cdots hP_s \, v^{(n)}.$$

As the action of $U(h \otimes A)$ commutes with the action of $S^n(A)$, this map is 
a homomorphism of algebras. Then we can describe it explicitly observing that
the generators $P \in A = T^1(A)$ are mapped to
${\rm Sym}(P \bigotimes 1 \bigotimes \dots \bigotimes 1)$.

Now consider the augmentation $\varepsilon_\szn$ of $S^n(A)$ defined by
$$\varepsilon_\szn(R) = R\left(\Z^{(1)}_1, \dots, \Z^{(1)}_n\right).$$
We can pull it back to the augmentation
$\epsilon_\szn$ of $T^*(A)$. Then for the generators $P \in A$ we have
\begin{equation}\label{aug}
\epsilon_\szn(P) = \sum_{j=1}^n P\left(\Z^{(1)}_j\right).
\end{equation}
As this is exactly the action of $hP \in h \otimes A$ on $v_\szn$, the augmentation ideal
$J_\szn$ acts on the image of $v^{(n)}$ in $W_M(\szn)$ by zero.

Now we have the surjective map $\phi: \natm / J_\szn \cdot \natm \to W_M(\szn)$.
The module in the left hand side is generated by the image $v'$ of $v^{(n)} \in \natm$.
Note that $v'$ satisfies the
same conditions~\eqref{hwl} 
as $v_\szn$, namely $e \otimes A$ acts on $v^{(n)}$ and therefore on $v'$ by zero,
and $h \otimes A$ acts on $v'$ by the augmentation~\eqref{aug}.
As $\phi(v') = v_\szn$, the universality property of
$W_M(\szn)$ implies that $\phi$ is an isomorphism.
\end{proof}

\begin{prp}\label{quo}
Any quotient $W$ of $W_A(n\omega)$ such that $\dim W^{\n\omega} =1$ is
indeed a
quotient of $W_M(\szn)$ for a certain $n\omega$-multiset $\szn$.
\end{prp}

\begin{proof}
Let $v$ be the image of $v_{n\omega}$ in $W$. We know that $h \otimes A$
acts on $v$ by a certain functional on $A$, so the tensor algebra 
$T^*(A) \cong U(h \otimes A)$ acts on $v$ by a certain augmentation $\epsilon$.

Note that the defined above map $T^*(A) \to S^n(A)$ is surjective and that its kernel is
contained
in the augmentation ideal of $\epsilon$.
Therefore $\epsilon$ is the pull back of
a certain augmentation $\varepsilon$ of $S^n(A)$. It means  that $S^n(A)$
acts on $v$ by $\varepsilon$. We have $\varepsilon \ne 0$ as the action of $1
\in S^n(A)$ on $v$ is not zero.

So the kernel of $\varepsilon$ is a maximal ideal in $S^n(A)$, i.e. a point on
$S^n(M)$, that is a set of $n$ points on $M$ (not necessarily distinct).
Denote
them by $\Z^{(1)}_j$, $j=1\dots n$. Then the vector $v$ satisfies the
same conditions~\eqref{hwl} as $v_\szn$. Universality property implies the proposition.
\end{proof}

\begin{thm}\label{univ}
Let $\g$ be an arbitrary semisimple Lie algebra. Then
any quotient $W$ of $W_A(\lambda)$ such that $\dim W^{\lambda} =1$ is
indeed a
quotient of $W_M(\sz)$ for a certain $\lambda$-multiset $\sz$.
\end{thm}

\begin{proof}
Let $v$ be the image of $v_\lambda$ in $W$. For $i=1\dots r$ by 
$sl_2^{(i)}$ denote the
subalgebra of $\g$ generated by $e_i$, $h_i$, $f_i$.

Let $W^{(i)}$ be the subspace of $W$ generated by $v$ under the action of
$sl_2^{(i)} \otimes A$. Then $W^{(i)}$ satisfies the conditions of
Proposition~\ref{quo}. So $h_i \otimes A$ acts on $v$
according to~\eqref{hwl}
for a certain set $\{\Z^{(i)}_j\}_{j=1\dots \lambda_i}$.
Combining
these sets for $i=1\dots r$ into a $\lambda$-multiset,
we obtain the statement of the theorem.
\end{proof}

\subsection{Beyond $sl_2$}

\def \szo {{\{\Z\}_{n\omega_1}}}

Let us first consider the case $\g = sl_{r+1}$.
By $e_{ij}$, $f_{ij}$, $1\le i<j\le r+1$ denote the
matrix units corresponding to positive and
negative roots respectively.

Let $V_{r+1}$ be the vector representation. We can
similarly define $L_N(V_{r+1}) = V_{r+1} \otimes A$ and consider the
symmetric powers.

\def \natmr  {S^n \left(L_A(V_{r+1})\right)}

\begin{prp}\label{gdecmp}
We have 
$$
\natmr \cong \bigoplus_{{i_0 + i_1+ \dots + i_{r} = n}\atop{i_j\ge 0}} 
\natmr^{(i_0 - i_1, i_1 - i_2, \dots, i_{r-1}-i_r)},\ \ \mbox{where}
$$
$$\natmr^{(i_0 - i_1, i_1 - i_2, \dots, i_{r-1}-i_r) } \cong
S^{i_0}(A) \bigotimes \dots \bigotimes S^{i_r}(A)$$
and $h_j \otimes 1$ acts on $\natmr^{(s_1, \dots, s_r)}$ by the scalar $s_j$. 
\end{prp}

\begin{proof}
Similar to Proposition~\ref{decmp}.
\end{proof}

\begin{thm}\label{giso}
We have the natural isomorphism of $sl_{r+1} \otimes A$-modules
$$W_A(n\omega_1) \cong \natmr.$$
\end{thm}

\begin{proof}
Similar to Theorem~\ref{iso}. The basis is formed by monomials
\begin{equation}\label{nmon}
(f_{12}P_{j^1_1} \cdots f_{12}P_{j^1_{i_1}}) \cdot \dots
\cdot
(f_{1,r+1}P_{j^r_1} \cdots f_{1,r+1}P_{j^r_{i_r}}) \cdot
(h_1 P_{j^0_1} \cdots h_1 P_{j^0_{i_0}}) \, v^{(n)},
\qquad j^m_1 \ge \dots \ge j^m_{i_m},
\end{equation}
with $i_0+ \dots + i_r = n$.

The analog of Proposition~\ref{bas} is straightforwardly similar. 
Repeating the proof of Proposition~\ref{span}, 
let us show that these monomials span the Weyl module.

By Lemma~\ref{pgr} we have the restriction $i_1 + \dots + i_r \le n$. Let us fix the grading,
that is fix the integers $i_1, \dots, i_r$ (but not $i_0$). 
Similarly introduce the partial order comparing $i_0$, and then
the number of  $j^r_k$ equal to one.

Then each monomial~\eqref{nmon}
 can be obtained as the highest monomial in
$$(e_{1,r+1} P_{j^0_1} \cdots e_{1,r+1}  P_{j^0_{i_0}})  \cdot f_{1,r+1}^{i_0} \cdot
(f_{12}P_{j^1_1} \cdots f_{12}P_{j^1_{i_1}}) \cdot \dots \cdot
(f_{1,r+1}P_{j^r_1} \cdots f_{1,r+1}P_{j^r_{i_r}}) \, v^{(n)}$$
with respect to this order.
So for $i_0 + i_1 + \dots + i_n >n$ 
each monomial~\eqref{nmon} is the highest monomial of a relation, 
therefore by Lemma~\ref{ordspan} (ii)
the other monomials~\eqref{nmon} span the module.
For $i_0 + i_1 + \dots + i_n <n$ we just set the additional indexes $j^0_s$ equal to one.
\end{proof}

Define the action of $S^n(A)$ and the ideal $J_\szo$ as above.

\begin{thm}\label{giso2}
For any $\szo$ we have the natural isomorphism of $sl_{r+1}
\otimes
A$-modules
$$W_M(\szo) \cong \natmr / J_{\szo}\cdot \natmr$$
\end{thm}

\begin{proof}
Straightforwardly similar to Theorem~\ref{iso2}.
\end{proof}

In the case of an arbitrary Lie algebra $\g$ 
the following analog of the
conjecture \cite{CP} about global Weyl modules looks reasonable.

\begin{cnj}
Let $V^{(i)}$, $i = 1 \dots r$ be the fundamental representations of
$\g$. By $v^{(i)} \in V^{(i)}$ denote the highest weight vectors.
Then the module $W_A(\lambda)$ over the algebra $\g \otimes A$ is
isomorphic
to
the submodule of
$$ \left(V^{(1)}\otimes A\right)^{\otimes \lambda_1} \bigotimes \dots
\bigotimes
\left(V^{(r)}\otimes A\right)^{\otimes \lambda_r}$$
generated by $\left(v^{(1)} \otimes 1\right)^{\otimes \lambda_1}
\bigotimes \dots \bigotimes \left(v^{(r)}
\otimes 1\right)^{\otimes \lambda_r}$.
\end{cnj}

We also expect some results similar to \cite{CK}.

\section{Case of polynomial algebra}

\def \bx {{\mathbf x}}
\def \ol {{\{0\}}}
\def \fr {{\mathcal F}_n^{r+1}}
\def \cnv {{/\cdot}}

\subsection{General case}

Recall that by Corollary~\ref{elw} any local Weyl module on a
non-singular variety of dimension $d$ is isomorphic to the tensor
product of some 
modules $W_{\CC^d}(\ol_{\lambda^i})$. Let us describe them for $\g =sl_{r+1}$ and
the weight equal to $n \omega_1$ in more detail.

Let $M = \CC^d$, then $A = \CC[x^1, \dots, x^d]$. We will write it as
$A = \CC[\bx]$,
where $\bx = (x^1, \dots, x^d) \in
\CC^d$. 
Also fix the notation $\CC[\bx_1, \dots, \bx_n]$ for the algebra of polynomials in
$n$ multivariables (that is in $dn$ usual variables) and the notation
$\CC[\bx_1, \dots, \bx_n]_+$ for the ideal of polynomials with zero
at the origin. 

Note that the symmetric group $\Sigma_n$ acts on $\CC[\bx_1, \dots, \bx_n]$ by 
permutation of multivariables. Then the space of invariants
$\CC[\bx_1, \dots, \bx_n]^{\Sigma_n}$ generalizes the notion of symmetric polynomials.

Introduce the {\em Frobenius transformation} $\fr$ 
from the Grothendieck ring of the symmetric group $\Sigma_n$ to the Grothendieck ring of $sl_{r+1}$
mapping any $\Sigma_n$-module $\pi$ to the $sl_{r+1}$-module 
$(V_{r+1}^{\otimes n} \otimes \pi)^{\Sigma_n}$. Here $sl_{r+1}$ acts on the 
factors $V_{r+1}$, and $\Sigma_n$
acts on $V_{r+1}^{\otimes n}$ by permuting the factors as well as on $\pi$. 

To be explicit, $\fr$ maps the irreducible representation of $\Sigma_n$ corresponding to a
partition $\xi$ to the irreducible representation of $sl_{r+1}$ with highest weight
$(\xi_1 - \xi_2, \dots, \xi_r-\xi_{r+1})$ if $\xi_{r+2}=0$ and to zero otherwise.

For an algebra $A$ and a subalgebra $B$ introduce the notation $A\cnv B = A/(B\cdot A)$.

\begin{crl}\label{frc}
 For $\g = sl_{r+1}$ we have an isomorphism of $sl_{r+1}$-modules
$$W_{\CC^d}(\ol_{n\omega_1}) \cong \fr \left(
\CC[\bx_1, \dots, \bx_n]\cnv \CC[\bx_1, \dots, \bx_n]_+^{\Sigma_n}\right).$$
\end{crl}

\begin{proof}
 Note that $A^{\otimes n} \cong \CC[\bx_1, \dots, \bx_n]$ 
and $S^n(A) \cong \CC[\bx_1, \dots, \bx_n]^{\Sigma_n}$. Therefore
$$\natmr \cong \left (V_{r+1}^{\otimes n} \otimes \CC[\bx_1, \dots, \bx_n]\right)^{\Sigma_n}$$ 
and $J_{\ol_{\n\omega_1}} \cong  \CC[\bx_1, \dots,
\bx_n]^{\Sigma_n}_+$. So by Theorem~\ref{giso2} the quotient $W_{\CC^d}(\ol_{n\omega_1})$
is exactly as proposed.
\end{proof}

Also let us calculate weights multiplicities of this module.

\begin{crl}\label{csym}
 For $\g = sl_{r+1}$ we have
$$W_{\CC^d}(\ol_{n\omega_1})^{(i_0 - i_1, i_1 - i_2, \dots, i_{r-1}-i_r) }
\cong \CC[\bx_1, \dots, \bx_n]^{\Sigma_{i_0}\times \dots \times \Sigma_{i_r}}\cnv 
\CC[\bx_1, \dots, \bx_n]_+^{\Sigma_n},$$
where $i_0 + \dots + i_r = n$ and $\Sigma_{i_0}\times \dots \times \Sigma_{i_r}$ 
acts as a subgroup of $\Sigma_n$.
\end{crl}

\begin{proof}
It follows in the similar way from  Theorem~\ref{giso2} and
Proposition~\ref{gdecmp}.
\end{proof}

So investigation of Weyl modules in this case is reduced to a question
about the multi-variable analog of symmetric functions.

\subsection{One-dimensional case}

This calculation is based on the following statement.

\begin{thm} {\em (Chevalley)}
The representation of $\Sigma_n$ in the quotient $\CC[x_1, \dots, x_n]\cnv
\CC[x_1, \dots, x_n]_+^{\Sigma_n}$ is isomorphic to the regular one.
\end{thm}

In \cite{CP} it is shown that in $sl_2$ case we have
$$\dim W_{\CC}(\ol_{n\omega})
= 2^n, \qquad \dim W_{\CC}(\ol_{n\omega})^{n-2i} = \bin{n}{i}.$$

Let us slightly generalize this statement.

\begin{prp} For $\g = sl_{r+1}$ we have
$$ \dim W_\CC(\ol_{n\omega_1}) = (r+1)^n, \qquad
\dim W_\CC(\ol_{n\omega_1})^{(i_0 - i_1, i_1 - i_2, \dots, i_{r-1}-i_r) }
=
\frac{n!}{i_0! i_1! \dots i_r!}$$
for all  $i_0 + \dots + i_r = n$.
\end{prp}

\begin{proof}
Let us prove the second formula, then the first one follows. By
Corollary~\ref{csym} we have

$$W_\CC(\ol_{n\omega_1})^{(i_0 - i_1, i_1 - i_2, \dots, i_{r-1}-i_r) }
\cong \CC[x_1, \dots, x_n]^{\Sigma_{i_0}\times \dots \times \Sigma_{i_r}}\cnv \CC[x_1, \dots, x_n]_+^{\Sigma_n}$$

So we have the
subspace of invariants in the regular representation
$\CC[\Sigma_n]$
with respect to the subgroup
$\Sigma_{i_0} \times \dots \times \Sigma_{i_r}$. Then the dimension is
equal to the index of the subgroup.
\end{proof}

\subsection{Two-dimensional case}

This calculation is based on the following result from \cite{H}.

Recall that a function $f : \{1, \dots, n\} \to \{1, \dots, n\}$ is called
a {\em parking} function if $|f^{-1}(\{1, \dots, k\})| \ge k$ for $k = 1\dots n$.

The group $\Sigma_n$ acts on the set $PF_n$ of parking functions by
permutation
on the domain. Let $\CC PF_n$ be the permutation representation of
$\Sigma_n$
on parking functions.

\begin{thm} \label{th}\ {\rm (See \cite{H}, Theorem~3.10 together with
Proposition~3.13)}
Let $\bx_i \in \CC^2$. Then
the representation of $\Sigma_n$ in the quotient $\CC[\bx_1, \dots,
\bx_n]\cnv
\CC[\bx_1, \dots, \bx_n]_+^{\Sigma_n}$ is isomorphic to the product $\CC
PF_n
\otimes \varepsilon$, where $\varepsilon$ is the sign representation of
$\Sigma_n$.
\end{thm}

\def \ca {{\mathcal A}}

To describe $\CC PF_n$ more explicitly introduce the class of
integer sequences
$$\ca_n = \{ (a_1, \dots, a_n) | \ a_i \ge 0, \ a_1+ \dots + a_n = n,\ a_1 +
\dots + a_k
\ge k\ \mbox{for}\ k=1\dots n\}.$$

\begin{prp}\label{indp} We have
$$\CC PF_n = \bigoplus_{\ca_n} \ind_{\Sigma_{a_1} \times \dots \times
\Sigma_{a_n}}^{\Sigma_n} \CC.$$
\end{prp}

\begin{proof}
For a sequence $(a_1, \dots a_n)$ with
$a_i \ge 0$, $a_1+ \dots + a_n = n$
consider the set of 
functions $f : \{1, \dots, n\} \to \{1, \dots, n\}$ 
such that $|f^{-1}(\{i\})| = a_i$, $i=1\dots n$. This set admits
the action of $\Sigma_n$ by permutation
on the domain and forms the permutation
representation isomorphic to
$$\CC[\Sigma_n / \Sigma_{a_1} \times \dots \times
\Sigma_{a_n}] \cong  \ind_{\Sigma_{a_1} \times \dots \times
\Sigma_{a_n}}^{\Sigma_n} \CC.$$
To complete the proof note that parking functions just correspond to 
$(a_1, \dots a_n) \in \ca_n$.
\end{proof}

\def \hs {{\mathcal R}}

By $\hs_n^s$ denote the set of subsets $H \subset \{1 \dots sn\}$ such
that
$$|H| = n, \qquad |H\cap \{1, \dots, k\}| \ge \frac{k}{s}.$$
Such subsets are usually called {\em Raney sequences} (see \cite{GKP},
\cite{S}).

\begin{thm}\label{mcalc} Let $\g = sl_{r+1}$.
\begin{enumerate}
\item
The dimension of the space
$W_{\CC^2}(\ol_{n\omega_1})$ is equal to the number of
elements in $\hs_{n+1}^{r+1}$.
\item
The dimension of the space
$W_{\CC^2}(\ol_{n\omega_1})^{(i_0 - i_1, i_1 - i_2, \dots, i_{r-1}-i_r) }$
is equal to the number of elements $H \in \hs_{n+1}^{r+1}$ such that
$$|\{ i \in H| \,i \equiv j+1 \, {\rm mod}\, (r+1)\}| = i_j,\ j=1\dots
r.$$
\end{enumerate}
\end{thm}

\begin{proof}
By Corollary~\ref{frc} and Theorem~\ref{th} we have
$$W_{\CC^2}(\ol_{n\omega_1}) \cong \fr \left(
\CC[\bx_1, \dots, \bx_n]\cnv \CC[\bx_1, \dots, \bx_n]_+^{\Sigma_n}\right) \cong
\fr (\CC PF_n\otimes \varepsilon),$$
so the Frobenius duality and Proposition~\ref{indp} imply  
\begin{eqnarray*}
W_{\CC^2}(\ol_{n\omega_1}) \cong
\left(V_{r+1}^{\otimes n}\otimes \CC PF_n
\otimes \varepsilon\right)^{\Sigma_n} &\cong& \hom_{\Sigma_n}
\left(\varepsilon
\otimes (V_{r+1}^*)^{\otimes n}, \CC PF_n\right) \cong \\
\cong
\bigoplus_{\ca_n}
\hom_{\Sigma_n}
\left(\varepsilon
\otimes (V_{r+1}^*)^{\otimes n},\ind_{\Sigma_{a_1} \times \dots \times
\Sigma_{a_n}}^{\Sigma_n} \CC \right) &\cong& \bigoplus_{\ca_n}
\hom_{\Sigma_{a_1} \times \dots \times
\Sigma_{a_n}}\left(\varepsilon
\otimes (V_{r+1}^*)^{\otimes n}, \CC \right)
\cong \bigoplus_{\ca_n} \bigotimes_{i=1}^n \left(\wedge^{a_i}
V_{r+1}\right).
\end{eqnarray*}

Let $v$ be the highest weight vector of $V_{r+1}$. Choose the graded basis
$v_j \in V_{r+1}$, $j=1\dots r+1$, by setting $v_j = f_{1,j+1}\,v$ for
$j\le r$ and $v_{r+1}=v$. Then the graded monomials spanning
$W_{\CC^2}(\ol_{n\omega_1})$ look like
$$(v_{j^1_1}\wedge \dots \wedge v_{j^1_{a_1}})\bigotimes \dots \bigotimes
(v_{j^n_1}\wedge \dots \wedge v_{j^n_{a_n}}), \qquad 1 \le j^s_t \le
r+1, \qquad (a_1,\dots, a_n) \in \ca_n.$$

With this element we associate the set
$$H = \{1, (r+1)(s-1) + j^s_t + 1| \,
t = 1\dots a_s, \  s=1\dots n \}.$$

Then the condition $H \in \hs^{r+1}_{n+1}$ is equivalent to
$(a_1,\dots, a_n)  \in
\ca_n$. So we have a bijection between $\hs^{r+1}_{n+1}$ and graded
monomials. Concerning the grading note that the number of factors $v_j$
for $j= 1\dots r$
in a graded monomial is equal to the number of elements in the
corresponding set $H$ equivalent to $j+1$ modulo $r+1$.
\end{proof}

\begin{crl}
For $\g = sl_2$ we have
$$\dim W_{\CC^2}(\ol_{n\omega}) = \frac{(2n+2)!}{(n+1)!\,
(n+2)!},$$
that is the Catalan numbers and
$$\dim W_{\CC^2}(\ol_{n\omega})^{n-2i} = \frac{(n+1)! \,n!}{(n-i+1)!\,
(n-i)!\,(i+1)! \,i!}, \qquad i = 0\dots n,$$
that is the Narayana numbers.
\end{crl}

\begin{proof}
Let $C_0(z) = z^{-1}$ and
$$C_{n+1}(z) = \sum_{i=0}^n z^{n-2i}\cdot \frac{(n+1)! \,n!}{(n-i+1)!\,
(n-i)!\,(i+1)! \,i!}.$$
Then we have (see \cite{S}, Section~6)
$$C_{n+1}(z) = \sum_{s=0}^n C_s(z^{-1}) \cdot C_{n-s}(z).$$

To prove the corollary it is enough to check this relation for
$${\rm ch\,} W_{\CC^2}(\ol_{n\omega}) = \sum_{i=0}^n z^{n-2i} \cdot
\dim W_{\CC^2}(\ol_{n\omega})^{n-2i}.$$
By Theorem~\ref{mcalc} we can work with Raney sequences, namely let
$R_n^{(i)}$ be the number of sets $H \in \hs^2_n$ such that $H$ has
exactly $i$ even elements. Let $R_{0}(z) = z^{-1}$ and
$$R_{n+1}(z) = \sum_{i=0}^n z^{n-2i}\cdot R_{n+1}^{(i)}.$$
Then it is enough to show that $R_{n+1}(z) = \sum_{s=0}^n R_s(z^{-1})
\cdot R_{n-s}(z)$.

Take $H \in \hs_{n+1}^2$.
Let $s$ be
the first non-negative integer such that
$|H\cap \{1, \dots, 2s+2\}| = s+1$. Then we can split $H$ into
$$H_1 = \{i \in 1\dots 2s |\, i+1 \in H\}, \qquad H_2 = \{i \in 1\dots
2n-2s |\, i+2s+2 \in H\}.$$
Clearly $H_1 \in \hs^2_s$ and $H_2 \in \hs^2_{n-s}$.
On the other hand, if
$H_1 \in \hs^2_{s}$ and $H_2 \in \hs^2_{n-s}$ then we can glue them
into the set
$$H = \{1\} \cup \{i+1 |\, i \in H_1\} \cup \{i+2s+2 |\, i \in H_2\} \in \hs_{n+1}^2$$
Concerning the grading note that
the set of even elements in $H$ is formed by odd elements in $H_1$ and
even elements in $H_2$.
And the inversion of the formal variable $z$ in the recurrence relation
means exactly that we count odd elements instead of even ones.
\end{proof}

\begin{crl}
For $\g = sl_{r+1}$ we have
$$\dim W_{\CC^2}(\ol_{n\omega_1}) =
\frac{((r+1)(n+1))!}{(n+1)!(r(n+1)+1)!}.$$
That is the higher Catalan numbers.
\end{crl}

\begin{proof}
Similar to $sl_2$ case. See \cite{GKP}, Section 7.5 concerning Raney
sequences.
\end{proof}

So the weight multiplicities of these modules 
can be considered as higher Narayana numbers.

\subsection{A conjecture for higher dimensions}

\begin{cnj} For $\g = sl_2$ we have
$$\dim W_{\CC^d}(\ol_{n\omega})^{n-2i} = \frac{(d+n-1)! \cdots
(d+n-i)!}{(n-1)! \cdots (n-i)!}\cdot \frac{2! \cdots
(i-1)!}{d! \cdots (d+i-1)!}, \qquad i=0\dots n,$$
that equal to the dimension of the 
irreducible representation of $sl_n$ with highest weight $d\omega_i$ for $i=1\dots n-1$  
(and equal to one for $i=0,\, n$).
\end{cnj}

Then the dimensions of Weyl modules $W_{\CC^d}(\ol_{n\omega})$ form so called
Hoggatt sequences (see \cite {FA}).

\section*{Appendix}

Here we perform a formal computation used in this paper. 
We take a commutative algebra $A$ with unit and consider the Lie algebra $sl_2 \otimes A$.
Let $P_1, \dots P_n \in A$, and suppose that the vector $v$ is such that 
$(e\otimes P_i) v = (h \otimes P_i)v =0$ for all $i$. 
Let us write $f$ instead of $f \otimes 1$. Our aim is to calculate
$$(e\otimes P_1) \dots (e\otimes P_n) f^m v.$$

\def \cg {{\mathcal G}}
\def \cm {{\mathcal S}}

Suppose we have $m$ indistinguishable glasses (set $\cg_m$) 
and $n$ distinguishable ingredients ($I_1$, \dots, $I_n$).
By {\em cocktail serving} denote a way $S$ to distribute
ingredients between glasses (each ingredient into only one glass). 
By $\cm(m,n)$ denote the set of all possible cocktail serving.

We now associate an element of $U(sl_2 \otimes A)$ with 
a set of elements $P_1$ \dots, $P_n \in A$ and a cocktail serving $S$.
First for each glass $G\in \cg_m$ we take $P(G) \in A$ equal to the product of
$P_i$ for all $I_i$ from the glass $G$. 
Then we set $f(S) = \prod_{G \in \cg_m} (f\otimes P(G))$.

\begin{prp}\label{martini}
We have
\begin{equation}\label{mart}
(e\otimes P_1) \dots (e\otimes P_n) f^{n+m} v = \sum_{S \in \cm(m,n)}
c(S) f(S) v,
\end{equation}
where $c(S)$ are integers and $(-1)^n c(S) >0$ for $m>0$.
\end{prp}

\begin{proof}
Induction on $n$. For $n=0$ we have empty glasses and the only constant equal to one.
Concerning the step, note that
\begin{eqnarray*}
(e\otimes P) (f\otimes Q_1) \dots (f\otimes Q_n) v = 
\sum_{i=1}^n (f\otimes Q_1) \dots (f\otimes Q_{i-1}) (h \otimes PQ_i) 
(f\otimes Q_{i+1}) \dots (f\otimes Q_n)v =\\
= - 2\sum_{i< j} 
(f \otimes PQ_iQ_j) \prod_{k \ne i, j} 
(f\otimes Q_k)v= 
- \sum_{i\ne j} 
(f \otimes PQ_iQ_j) \prod_{k \ne i, j} 
(f\otimes Q_k)v.
\end{eqnarray*}

\def \mix {{\rm Mix}}

For $G_1,\,G_2 \in \cg_{m}$ introduce the {\em mixing} map
$\mix_{n+1}^{G_1,G_2} : \cm(m,n) \to \cm(m-1,n+1)$ that identifies the glass $G_1$ 
with the glass $G_2$,
puts there all their 
ingredients and finally adds there the ingredient $I_{n+1}$. Then for any
$S \in \cm(m,n-1)$ we have
$$(e\otimes P_n) f(S) v = - \sum_{G_1\ne G_2\in \cg_m} f(\mix_n^{G_1,G_2} (S)) v.$$

By induction we have~\eqref{mart}, but it remains to show that $(-1)^n c(S) >0$.
Note that $c(S)$ is the sum of $-c(S')$ such that $S' \in \cm(m+1,n-1)$
and $S$ is obtained from $S'$ by mixing. Therefore  $(-1)^n c(S) \ge 0$
and we have to show that the mixing map is surjective.
Consider $S'$ obtained from $S$ as follows.
We remove the glass $G$ containing $I_n$ and instead take the 
glasses $G_1$ and $G_2$ containing the other (excluding $I_n$) ingredients from $G$. Then 
we have $S = \mix_n^{G_1,G_2} (S')$.
\end{proof}

\end{document}